\documentclass[12pt]{article}
\usepackage{amssymb,graphicx}

\newcommand{\edo}{\end{document}}
\newcommand{\extraref}[1]{}  

\newcommand{\R}{{\mathbb R}}  

\def\ben{\begin{enumerate}}
\def\een{\end{enumerate}}

\newtheorem{theorem}{Theorem}
\newtheorem{itlemma}{Lemma}[section] 
\newtheorem{itproposition}[itlemma]{Proposition}
\newtheorem{itcorollary}[itlemma]{Corollary}
\newtheorem{itremark}[itlemma]{Remark}
\newtheorem{itdefinition}[itlemma]{Definition}
\newtheorem{itexample}[itlemma]{Example}

\newenvironment{lemma}{\begin{itlemma}\rm}{\end{itlemma}} 
\newenvironment{remark}{\begin{itremark}\rm}{\end{itremark}} 
\newenvironment{corollary}{\begin{itcorollary}\rm}{\end{itcorollary}}
\newenvironment{proposition}{\begin{itproposition}\rm}{\end{itproposition}}
\newenvironment{definition}{\begin{itdefinition}\rm}{\end{itdefinition}}
\newenvironment{example}{\begin{itexample}\rm}{\end{itexample}}
\newenvironment{proof}{\noindent {\em Proof}.\
}{\hspace*{\fill}$\halmos$\medskip}

\newcommand{\text}[1]{\hbox{\rm \ #1\ \/}}
\newcommand{\be}[1]{\begin{equation}\label{#1}}
\newcommand{\ee}{\end{equation}}
\newcommand{\beqn}{\begin{eqnarray*}}
\newcommand{\eeqn}{\end{eqnarray*}}
\newcommand{\beq}{\begin{eqnarray}}
\newcommand{\eeq}{\end{eqnarray}}
\newcommand{\bl}[1]{\begin{lemma}\label{#1}}
\newcommand{\ble}[1]{\begin{lemmaex}\label{#1}}
\newcommand{\br}[1]{\begin{remark}\label{#1}}
\newcommand{\bt}[1]{\begin{theorem}\label{#1}}
\newcommand{\bd}[1]{\begin{definition}\label{#1}}
\newcommand{\bp}[1]{\begin{proposition}\label{#1}}
\newcommand{\bc}[1]{\begin{corollary}\label{#1}}
\newcommand{\bfact}[1]{\begin{fact}\label{#1}}
\newcommand{\ber}[1]{\begin{exercise}\label{#1}}
\newcommand{\bex}[1]{\begin{example}\label{#1}}
\newcommand{\bem}[1]{\begin{example}\label{#1}}  
\newcommand{\ec}{\mybox\end{corollary}}
\newcommand{\efact}{\mybox\end{fact}}
\newcommand{\eer}{\mybox\end{exercise}}
\newcommand{\eex}{\mybox\end{example}}
\newcommand{\eem}{\mybox\end{example}}
\newcommand{\el}{\mybox\end{lemma}}
\newcommand{\ele}{\mybox\end{lemmaex}}
\newcommand{\er}{\mybox\end{remark}}
\newcommand{\et}{\qed\end{theorem}}
\newcommand{\ed}{\mybox\end{definition}}
\newcommand{\ep}{\mybox\end{proposition}}
\newcommand{\epr}{\end{proof}}
\newcommand{\bpr}{\begin{proof}}

\newcommand{\ecs}{\end{corollary}}
\newcommand{\eers}{\end{exercise}}
\newcommand{\eexs}{\end{example}}
\newcommand{\eems}{\end{example}}
\newcommand{\els}{\end{lemma}}
\newcommand{\eles}{\end{lemmaex}}
\newcommand{\ers}{\end{remark}}
\newcommand{\ets}{\end{theorem}}
\newcommand{\eds}{\end{definition}}
\newcommand{\eps}{\end{proposition}}
\newcommand{\halmos}{\rule{1ex}{1.4ex}}
\newcommand{\qed}{\hfill \halmos} 
\newcommand{\mybox}{\hfill $\Box$} 

\newcommand{\st}{\; | \;}

\newcommand{\edeq}{\vskip-0.7cm\mybox\end{definition}}

\def\bi{\begin{itemize}}
\def\ei{\end{itemize}}

\newcommand{\interK}{{\mbox{int}(K) }}

\newcommand{\Ker}{{\mbox{ker}\,}}

\newcommand{\Su}{P}
\newcommand{\Pd}{Q}
\newcommand{\ESu}{C}
\newcommand{\FPr}{D}

\begin{document}
\title{A Global Convergence Result for Strongly Monotone Systems with Positive
  Translation Invariance} 
\author{David Angeli\\
Dip. di Sistemi e Informatica\\
Universit\'a di Firenze\\
Via di S. Marta 3, 50139 Firenze, Italy\\
email: 
{\tt angeli@dsi.unifi.it}
\and Eduardo D. Sontag\footnote{corresponding author;
Phone +1.732.445.3072; FAX +1.206.338.2736}\\
Department of Mathematics\\
Rutgers University\\
Piscataway, NJ 08854-8019, USA\\
{\tt http:/\null\hskip-2pt/www.math.rutgers.edu/$\,{}_{\textstyle\tilde{}}\;${sontag}}\\
email:{\tt sontag@control.rutgers.edu}}
\maketitle

\textbf{Abstract} 

We show that strongly monotone systems of ordinary differential equations
which have a certain translation-invariance property are so that all
solutions converge to a unique equilibrium.
The result may be seen as a dual of a well-known theorem of Mierczynski for
systems that satisfy a conservation law.
An application to a reaction of interest in biochemistry is provided as an
illustration.


\section{Introduction and Motivations}

We recall that a dynamical system is said to be \emph{monotone} whenever
its state space $X$ is endowed with a partial order $\succeq$ and
the forward flow preserves the order.  
In other words,
for each ordered pair of initial conditions $\xi _1 \succeq \xi _2$,
solutions remain ordered:
$\varphi_t (\xi _1) \succeq \varphi_t (\xi _2 )$ for all $t\geq 0$.
See~\cite{Smith} for a discussion and many basic theorems, as well
as the recent excellent exposition~\cite{Hirsch-Smith}.
A special and most interesting case is when the partial order is induced by a
positivity cone, i.e.\ a closed subset
$K$ of a Banach space $B$ containing $X$ such that
$K + K \subset K$,
$K \subset \alpha K$ for all $\alpha \geq 0$,
and $K \cup -K = \{0 \}$.
In this case, one defines a partial order by the rule that
$\xi _1 \succeq \xi _2$ whenever $\xi _1 - \xi _2 \in K$. 
Strict versions of the order are also possible, and particularly useful
whenever $K$ has non-empty interior: one defines
$\xi _1 \succ \xi _2$  if $\xi _1 \succeq \xi _2$ and 
$\xi _1 \neq \xi _2$, and the following even stronger notion: $\xi _1 \gg \xi _2$ if
$\xi _1 - \xi _2 \in  \interK$.
A \emph{strongly} monotone system is one for which the following holds:
\be{defstrong}
\xi _1 \succ \xi _2 \; \Rightarrow \; \varphi_t ( \xi _1 ) \gg \varphi_t (\xi _2) \qquad \forall \, t>0,
 \quad \forall \, \xi _1,\xi _2 \in X.
\ee
A key foundational result is Hirsch's Generic Convergence Theorem
(\cite{Hirsch,Hirsch2,hirsch,Hirsch-Smith,Smith}), which guarantees that, if 
solutions of such systems are bounded, then, generically, they converge to
the set of equilibria. 
Roughly speaking, more complex asymptotic behaviors are possible,
but are (if they exist at all) confined to a zero-measure set of initial
conditions.

Remarkably, under suitable additional assumptions, generic convergence to
equilibria can be made global, as is the case if, for instance, the
equilibrium is unique \cite{Smith},  sometimes not requiring strong
monotonicity \cite{JiFa,Dancer}, if the system is cooperative and
tridiagonal \cite{Smillie} or if, there exists a positive first-integral for
the system, as shown in Mierczinski's paper \cite{mier}.
Our main result may be viewed as a dual of the latter result,
and applies to strongly monotone systems which have the property of
translation invariance with respect to a positive vector.
Equilibria of such systems are never unique.
The result is roughly as follows.
For systems evolving on Euclidean spaces $\R^n$, we
will assume that for some $v \in  \interK$, and for all
$\lambda \in \R$, the following is true: 
\be{invariance}
 \varphi_t (\xi + \lambda v) = \varphi_t (\xi ) + \lambda v 
 \ee   
for all $t \in \R$ and all $\xi \in X$. 
Under strong monotonicity, we show that convergence to equilibria is global
for a suitable projection of the system.
We also show that for competitive systems, i.e.\ systems that are
strongly monotone under time-reversal, the same result holds.
Statements and proofs are in Section~\ref{secmain}.

We were originally motivated by proving a global convergence result for
certain chemical reaction systems which are not necessarily monotone.
There has been much interest in recent years in establishing such global
results, see for 
instance~\cite{feinberg,volpert,kunze,gouze,sysbio,LSUmonotone,craciun}.
In Section~\ref{secchem}, we show how to associate, to any chemical reaction
system, a new system of differential equations, evolving on a different
space (of ``reaction coordinates'') for which our techniques may sometimes be
applied, and we illustrate with a system of interest in biochemistry.

In the last section, we make some remarks on extensions and comment on the
duality with Mierczinski's theorem.

\section{Main Result}
\label{secmain}

We consider nonlinear dynamical systems of the following form:
\be{thesystem}
\dot {x} = f(x)
\ee
with states $x \in X \subset \R^n$, for some closed set $X$ which is the closure of its interior, and some locally
Lipschitz vector field $f: X \rightarrow \R^n$.
For each initial condition $\xi \in X$, we denote by $\varphi_t ( \xi )$ the corresponding
solution,
and we assume that $\varphi_t (\xi )$ is uniquely maximally defined (as an element
of $X$) for $t \in  I_{\xi }$, where $I_{\xi }$ is an interval in $\R$ which
contains $[0,+\infty )$ in its interior.
(In other words, the system is assumed to be forward
--but not necessarily backward-- complete.)

Furthermore, a closed cone $K \subset \R^n$ is given, with non-empty interior,
and the corresponding non-strict and strict partial orders are considered:
$\succeq, \succ, \gg$. In particular, we assume that (\ref{thesystem}) 
is \emph{strongly monotone} as in (\ref{defstrong}) and that solutions enjoy
the translation invariance property (\ref{invariance}) 
for some $v \in  \interK $,
which we take, without loss of generality, to have norm one.

Because of property (\ref{invariance}) it is natural to assume, and we will do
so, that the state space is invariant with respect to translation by $v$,
namely: 
\be{statespace}
x \in X \; \Rightarrow \; x + \lambda v \in X \quad \forall \, \lambda \in \R\,.
\ee

In order to state our main result, we require an additional definition.
Given any unit vector $v$, we introduce the linear mapping:
\[
\pi _v: \R^n\rightarrow \R^n : x\mapsto  x - (v' x)v
\]
(prime indicates transpose),
which amounts to subtracting the component along the vector $v$, that is,
an orthogonal projection onto $v^\perp$.  Since $(v' x)v=(vv')x$,
we can also write $\pi _vx = (I-vv')x$.
Note that $\pi _vv=0$.

\noindent\textbf{Definition}
Let $\xi \in X$ be given and consider the corresponding solution $\varphi_t (\xi )$. 
We say that $\varphi_t(\xi )$ \emph{is bounded modulo $v$} if 
$\pi _v( \varphi_t (\xi ) )$ is bounded as a function of $t$, for $t\geq 0$.

Notice that we are not asking for precompactness of $\varphi_t (\xi )$ (which, in
examples, will typically fail), but only of its projection.

\noindent\textbf{Remark}
Equivalently, the solution $\varphi_t(\xi )$ is bounded modulo $v$ if and only if
there exists some scalar function
$\beta (\xi ,t): X \times [0,+ \infty) \rightarrow \R$ such that $\varphi_t (\xi ) - \beta (\xi ,t) v$ is
bounded as a function of time $t$.
(Recall that $X$ is invariant under translations by $v$, so this difference is
again an element of $X$.)
One direction is clear, using $\beta (\xi ,t)=v'\varphi_t(\xi )$.  Conversely, suppose
that there is any such $\beta $.  Then:
$v'(\beta (\xi ,t)v)=\beta (\xi ,t)v'v = \beta (\xi ,t)|v|^2 = \beta (\xi ,t)$,
so
$\pi _v(\beta (\xi ,t)v) = \beta (\xi ,t) - (v'(\beta (\xi ,t)v)v = 0$,
Since $X$ is closed, the assumption is that the closure of
$\{\varphi_t(\xi )-\beta (\xi ,t)v, t\geq 0\}$ is compact.
Thus, since $\pi _v$ is continuous, the same holds for
$\pi _v(\varphi_t(\xi )) = \pi _v(\varphi_t(\xi )-\beta (\xi ,t)v)$.

We are now ready to state our main result. 

\bt{main}
Consider a forward complete nonlinear system, strongly monotone on
$X$.
Let
(\ref{thesystem}) enjoy positive translation invariance as in
(\ref{invariance}) with respect to some vector $v \in  \interK$,
and so that the state space $X$ is closed and invariant with respect to
translation by $v$ as in (\ref{statespace}). Then, every solution which is
bounded modulo $v$ 
is such that $\pi _v( \varphi_t (\xi ) )$ converges to an equilibrium.
Moreover, there is a unique such equilibrium.
\ets

Before addressing the technical steps of the proof, it is useful to provide
an infinitesimal characterization of translation invariance.
This is a routine exercise, but we include a proof for ease of reference.

\bl{infinitesimal}
A system (\ref{thesystem}) enjoys the translation invariance property
(\ref{invariance}) with respect to $v\in \R^n$ if and only if:
\be{invariancevectorfield}
x_1,x_2\in X, \; 
x_1 - x_2 \in \textrm{span} \{ v \} \; \Rightarrow \; f(x_1)=f(x_2) \,.
\ee
\els

Notice that, for differentiable $f$, yet another characterization is
that $v\in \Ker f_*(x)$ (Jacobian) at all states $x$.

\bpr
If the system is translation-invariant by $v$, and $x_2=x_1+\lambda v$, then
$\varphi_t(x_2)-\varphi_t(x_1)=\lambda v$.
Taking $(d/dt)|_{t=0}$, we obtain $f(x_1)=f(x_2)$.
We now show the sufficiency of the condition.
More generally, suppose that $L$ is a linear subspace of $\R^n$ such that
$x_1-x_2\in L \Rightarrow  f(x_1)=f(x_2)$; we will prove that
$\varphi_t(x_2)-\varphi_t(x_1)$ is constant if $x_1-x_2\in L$.

We first change coordinates with a linear map $T$ in such a manner that $L$
gets transformed into the span of the first $\ell=\mbox{dim}\,L$ 
canonical vectors $\tilde L=\{e_1,\ldots ,e_\ell\}$.
The transformed equations are
$\dot {\tilde x} = \tilde f (\tilde x)$, where
$\tilde f(\tilde x) = T f(T^{-1}\tilde x)$ and $\tilde x = Tx$.
We partition the state as $\tilde x = (y',z')'$, with $y$ of size $\ell$, and
write the transformed equations in block form:
\beqn
\dot y &=& \tilde f_1(y,z)\\
\dot z &=& \tilde f_2(y,z) \,.
\eeqn
Suppose that two vectors $\tilde x_1$ and $\tilde x_2$ are such that
$z_1=z_2$.  This means that $\tilde x_1-\tilde x_2\in \tilde L$.
Then, letting $x_i:=T^{-1}\tilde x_i$, we have that $x_1-x_2\in L$,
and therefore $f(x_1)=f(x_2)$ by assumption.
Thus also 
$\tilde f(\tilde x_1) = Tf(x_1)=Tf(x_2)=\tilde f(\tilde x_2)$.
In other words, $\tilde f$ is independent of $y$, and the transformed
equations in block form read:
\beqn
\dot y &=& \tilde f_1(z)\\
\dot z &=& \tilde f_2(z) \,.
\eeqn
Now pick any $x_1,x_2\in X$ such that $x_1-x_2\in L$.
Then, $\tilde x_1-\tilde x_2\in \tilde L$, i.e., $z_1=z_2$.
Let $y_i(t)$ and $z_i(t)$ denote the components of the solution of the
transformed differential equation with respective initial conditions
$\tilde x_i$, $i=1,2$.
Then, $z_1(t)=z_2(t)$ for all $t\geq 0$ (same initial conditions for
the second block of variables), which implies that $\dot y_1(t)=\dot y_2(t)$
for all $t$.
Therefore also $\dot {\tilde x}_1(t)=\dot {\tilde x}_2(t)$ for all $t$,
and back in the original coordinates we have that
$(d/dt)(\varphi_t(x_2)-\varphi_t(x_1))=0$, as desired.
\epr

In order to carry out the proof we first need the following Lemma.

\bl{Lyapunovstrict}
Let $v \in \interK$ be given, such that $|v|=1$. Then, the function:
\[ V(x) := \inf \{ \alpha \in \R: x \preceq \alpha v \} \]
is well defined and Lipschitz for $x \in \R^n$.
\els

\bpr We show first that for all $x$ there exists $\alpha $ such that 
$\alpha v \succeq x$.
We may equivalently check that $v \succeq x/\alpha $ for some $\alpha \not= 0$.
Since $x/ \alpha \rightarrow 0$ as $\alpha \rightarrow + \infty$, we may conclude that this is the case,
since, as is well known, 
$(-v,v):= \{ x: v \gg x \gg -v \}$
is an open neighborhood of the origin, for all $v \gg 0$
(in other words the topology induced by a positivity cone with non-empty
interior is equivalent to the standard topology in $\R^n$).
On the other hand, $\alpha v \prec x$ for all sufficiently small $\alpha $
(as $\alpha \rightarrow -\infty $, $(-x)/(-\alpha )\rightarrow 0$, so $(-x)/(-\alpha )\prec v$, that is,
$-x\prec-\alpha v$).
Therefore, $V(x)$ is well defined.
Moreover, since $K$ is closed and the feasible set of
$\alpha $'s is bounded from below,
the infimum is achieved and is actually a minimum, which
implies that $V(x)$ is a continuous function.
We can prove, moreover, that $V$ is Lipschitz, as follows.

We first pick an $\varepsilon >0$ such that $\varepsilon z \prec v$ for all unit vectors
$z$.  (Such an $\varepsilon $ exists, because $\varepsilon z\rightarrow 0$ uniformly on the unit sphere,
and $(-v,v)$ is a neighborhood of zero.)
Therefore, for each two vectors $x\not= y$, applying this observation to
$z = \frac{1}{|x-y|}(x-y)$, we have that
$x-y  \preceq  k|x-y| v$, where $k := 1/\varepsilon $, and the same holds if $x=y$.
Now, given any two $x,y$, we write
\[
x \;=\; x-y + y  \;\preceq\; k|x-y|v + V(y)v \;=\; ( k|x-y| + V(y) ) v
\]
which means that $V(x) \leq  k|x-y| + V(y)$, and therefore
$V(x) - V(y) \leq  k|x-y|$.
Since $x$ and $y$ were arbitrary, this proves that $V$ is Lipschitz
with constant $k$.
\epr

The next Lemma is crucial for proving our main result.

\bl{isLyapunov}
Let $\xi _1$ and $\xi _2$ in $X$ 
be arbitrary, and $V$ be defined according to
the previous Lemma \ref{Lyapunovstrict}. Then, for all $t >0$ it holds that:
\be{decrease}
V ( \varphi_t (\xi _1) - \varphi_t ( \xi _2 ) ) \; \leq \; V ( \xi _1 - \xi _2) \,,
\ee
and the inequality is strict whenever $\xi _1 - \xi _2  \notin \textrm{span} \{ v \}$.
\els

\bpr 
Let $\xi _1$ and $\xi _2$ be arbitrary. By definition of $V$, we have:
$\xi _1 \preceq \xi _2 + V( \xi _1 - \xi _2 ) v$.
By translation invariance and monotonicity then:
$\varphi_t (\xi _1) \preceq \varphi_t (\xi _2 + V( \xi _1, \xi _2) v ) =
 \varphi_t (\xi _2 ) + V( \xi _1, \xi _2 ) v$. It follows that
 $V( \varphi_t (\xi _1) - \varphi_t(\xi _2) ) \leq V( \xi _1 - \xi _2)$, as claimed.
In particular, whenever $\xi _1 - \xi _2 \notin \textrm{span} \{ v \}$ we have
 $\xi _1 \prec \xi _2 + V( \xi _1 - \xi _2 ) v$ and therefore, exploiting strong
monotonicity: 
 $\varphi_t (\xi _1) \ll \varphi_t ( \xi _2 + V(\xi _1 - \xi _2 ) v ) =
 \varphi_t (\xi _2 ) + V (\xi _1 - \xi _2) v$. In particular, then
 $V (\varphi_t (\xi _1) - \varphi_t( \xi _2 ) ) < V (\xi _1 - \xi _2 )$.
 \epr
 
Notice that, by the semigroup property for flows, Lemma~\ref{isLyapunov}
implies that the function
$t\mapsto V (\varphi_t(\xi _1)-\varphi_t(\xi _2))$ is nondecreasing.

We also prove a result for systems that are strongly monotone in reversed
time, meaning that
for every  
pair $\xi _1, \xi _2$ and every time $t<0$ such that $\varphi_t ( \xi _1 )$ and 
$\varphi_t (\xi _2)$ are well-defined the following implication holds:
\[ 
\xi _1 \succ \xi _2 \; \Rightarrow \; \varphi_t ( \xi _1 ) \gg \varphi_t ( \xi _2 ). 
\] 

\bc{isLyapunovcomp}
Let $\xi _1$ and $\xi _2$ in $X$ 
be arbitrary, and $V$ be defined according to
the previous Lemma \ref{Lyapunovstrict}. 
Assume that system (\ref{thesystem}) be
forward-complete, strongly monotone in reversed time over $X$,
and translation invariant
with respect to some $v \in \interK$; then, for all $t >0$ it holds that:
\be{increase}
V ( \varphi_t (\xi _1) - \varphi_t ( \xi _2 ) ) \; \geq \; V ( \xi _1 - \xi _2) \,,
\ee
and the inequality is strict whenever $\xi _1 - \xi _2  \notin \textrm{span} \{ v \}$.
\ec

We are now ready to prove the main result.

\emph{Proof of Main Result}

Let $\xi \in X$
be such that $\varphi_t ( \xi )$ is bounded modulo $v$. 
That is,
$\tilde{x}(t)  := \pi _v(\varphi_t(\xi )) = (I-vv')\varphi_t(\xi )$ 
is a bounded function of $t$.  
Notice that $\tilde{x}$ satisfies the following differential equation:
\be{newsystem}
 \dot { \tilde{x} } = (I - v v') f ( \varphi_t ( \xi ) ) = (I - v v' ) f ( \tilde{x} ) 
 \ee
where the last equality follows by translation invariance.
This is a new dynamical system,
with state space $\tilde{X}:= \{ \tilde{x} \in v^{\perp}:
 \exists \, \lambda \in \R: \tilde{x} + \lambda v \in X \}$, 
viz.\ the projection along $v$
of $X$ onto the vector-space $v^{\perp}$, and we will denote by $\tilde{\varphi}_t$
the corresponding flow. 
Notice that $\pi $ (we omit the subscript $v$ from now on), 
$\varphi_t$ and $\tilde{\varphi}_t$ are related in the following sense: 
\[
\pi \circ \varphi_t \;=\; \tilde{\varphi} _t \circ \pi  \,.
\]
Moreover, by translation invariance of $X$, we have 
$ \tilde{X} = X \cap v^{ \perp}$ and $X = \tilde{X} \oplus \textrm{span} \{ v \}$.

By the above considerations, it makes sense to speak about the $\omega $-limit set
$\omega (\tilde{x})$ of solutions of (\ref{newsystem}), which by the boundedness
assumption, will be a compact, non-empty invariant set.
We would like to show that $\omega ( \tilde{x} )$ is a single equilibrium.

We show uniqueness first.
An equilibrium $\widetilde x$ of~(\ref{newsystem}) satisfies that $f(\widetilde x)$ belongs to
the span of $v$, let us say $f(\widetilde x)=rv$.
Therefore, the function $z(t)=\widetilde x+ tf(\widetilde x)$ is a solution of the system
$\dot x=f(x)$, since its derivative satisfies:
\[
\dot  z(t) =  
f(\widetilde x)  = f(\widetilde x+(rt)v) = f(z(t)),
\]
where the second inequality is by (the infinitesimal characterization of)
translation invariance.
Since $z(0)=\widetilde x$, we have that $\varphi_t(\widetilde x)=\widetilde x+ tf(\widetilde x)$ for all $t$.
Assuming that $\tilde{x}_1$ and $\tilde{x}_2$ are two distinct equilibria
for (\ref{newsystem}),
we have that
$\varphi_t( \tilde{x}_i ) = \tilde{x}_i + t f( \tilde{x}_i)$ (for $i=1,2$).
Hence, for all $t > 0$:
\begin{eqnarray}
\label{firstineq}
 V( \tilde{x}_1 - \tilde{x}_2 ) &>& V ( \varphi_t ( \tilde{x}_1 ) - \varphi_t ( \tilde{x}_2 ) ) \nonumber \\
 &=& V \big ( \tilde{x}_1 - \tilde{x}_2  + [ f ( \tilde{x}_1 ) - f( \tilde{x}_2 ) ] t \big )
\end{eqnarray}
By a symmetric argument, however, 
\be{symmetriccase}
 V( \tilde{x}_2 - \tilde{x}_1 ) > V \big ( \tilde{x}_2 - \tilde{x}_1  + [ f ( \tilde{x}_2 ) - f( \tilde{x}_1 ) ] t \big )
\ee
which should hold again for all $t>0$. It is straightforward, from definition
of $V(x)$, that the function be increasing with respect to (positive)
translations along $v$. Hence, the inequality in (\ref{firstineq}) implies 
$f(\tilde{x}_1)-f(\tilde{x}_2) \prec 0$, 
while, the second inequality gives  
$f(\tilde{x}_1)-f(\tilde{x}_2 )\succ 0$. But this is clearly a contradiction. 

Let $\tau >0$ be arbitrary; consider the solutions of (\ref{thesystem})
corresponding to $\xi $ and $\varphi_\tau (\xi )$.
We claim that $\varphi_t ( \varphi_{\tau } (\xi ) ) - \varphi_t (\xi )$ is bounded.
In fact, denoting by $\tilde{\varphi}_t$ the corresponding projections onto
$\tilde{X}$ and exploiting Lemma \ref{infinitesimal}, we obtain:
\be{boundedness}
\varphi_t ( \varphi_{\tau } ( \xi ) ) - \varphi_t (\xi ) =
\varphi_{\tau } ( \varphi_t (\xi ) ) - \varphi_t (\xi ) = \int_t^{t+ \tau } f ( \varphi_s ( \xi ) ) \, ds
= \int_t^{t+ \tau } f ( \tilde{ \varphi}_s (\pi (\xi )) ) \, ds
\ee
and so $|\varphi_t(\varphi_\tau (\xi ))-\varphi_t(\xi )|\leq \tau M$. where $M$ is an upper bound on
the magnitude of the vector field $f$ on a compact set that contains
the trajectory $\pi _v(\varphi_t(\xi ))$.

Hence, $V ( \varphi_t ( \varphi_{\tau } (\xi ) ) - \varphi_t ( \xi ) )$ is lower-bounded, and, by
virtue of Lemma \ref{isLyapunov}, is decreasing. 
Therefore, it admits a limit $\bar{V}> - \infty$ as $t \rightarrow + \infty$.
We claim that
\be{aa}
\varphi_t ( \varphi_\tau (\xi )) - \varphi_t (\xi ) \rightarrow \textrm{span} \{ v \} \,.
\ee

Suppose that this claim is false.
Then, since, as we just proved, $\varphi_t(\varphi_\tau (\xi ))-\varphi_t(\xi )$ is bounded,
there will be a sequence of times $t_n\rightarrow \infty $ and an 
$\delta _0\not\in \textrm{span} \{ v \}$
such that  $\varphi_{t_n}(\varphi_\tau (\xi ))-\varphi_{t_n}(\xi )\rightarrow \delta _0$.
Moreover, by precompactness of $\pi \circ \varphi_t (\xi )$, we can pick a subsequence of
$\{t_n\}$, which we denote without loss of generality in the same way,
such that $\pi \circ \varphi_{t_n}(\xi )\rightarrow \tilde{x}_0$, for some vector $\tilde x_0$.

So the pair $[\tilde x_0,\delta _0]$ belongs to the following set:
\be{bigomega}
\Omega = \big \{ [ \tilde{x}, \delta ] : \exists \, t_n \rightarrow + \infty:
\pi \circ \varphi_{t_n} (\xi ) \rightarrow \tilde{x} \textrm{ and } \varphi_{t_n} ( \varphi_\tau ( \xi ) ) -
\varphi_{t_n} ( \xi ) \rightarrow \delta \big \} \,.
\ee
We show next that $\Omega $ satisfies the following invariance property:
\be{biginvariance}
 \forall \, [ \tilde{x}, \delta ] \in \Omega ,
\; \forall \, t \geq 0, \quad [ \tilde{ \varphi}_t ( \tilde{x} ) , \varphi_t ( \tilde{x} + \delta ) - 
\varphi_t ( \tilde{x} ) ] \in \Omega . 
\ee

Pick any $[ \tilde{x}, \delta ] \in \Omega $ and some sequence $\{t_n\}$ as in the
definition of $\Omega $, as well as any fixed $t>0$.
{}From $\tilde{x} = \lim_{n \rightarrow + \infty} \pi \circ \varphi_{t_n} ( \xi ) $ and
continuity of the flow, we have:
\be{firstremark}
\tilde{ \varphi}_t ( \tilde{x} ) = \lim_{n \rightarrow + \infty} \tilde{\varphi}_t ( \pi \circ \varphi_{t_n} ( \xi )  )
= \lim_{n \rightarrow + \infty} \pi \circ \varphi_{t+ t_n} ( \xi ) .
\ee
Moreover, 
\beqn
\delta &=& \lim_{n \rightarrow + \infty} \varphi_{t_n} ( \varphi_{ \tau } ( \xi )) - \varphi_{t_n} ( \xi )
\\ 
&=& \lim_{n \rightarrow + \infty} \varphi_{\tau } ( \varphi_{ t_n} ( \xi )) - \varphi_{t_n} ( \xi )  \\
&=&
\lim_{n \rightarrow + \infty} \varphi_{ \tau } \big ( \tilde{ \varphi }_{t_n} ( \pi ( \xi ) ) + 
[ v' \varphi_{t_n} ( \xi ) ] v \big ) - \tilde{\varphi}_{t_n} ( \pi (\xi ) ) - [v' \varphi_{t_n} ( \xi )] v
\eeqn
where the last equality follows 
from $\tilde \varphi_t(\pi (\xi ))=\pi (\varphi_t(\xi ))=\varphi_t(\xi )-[v'\varphi_t(\xi )]v$
applied to $t=t_n$.
Finally, from the equality $\varphi_\tau (\zeta +\lambda v)=\varphi_\tau (\zeta )+\lambda v$ applied to
$\zeta =\tilde\varphi_{t_n}(\pi (\xi ))$ and $\lambda =v'\varphi_{t_n}(\xi )$, this last expression 
gives that
\be{deltaprop}
\delta  \;=\;
\lim_{n \rightarrow + \infty} \varphi_{ \tau } \big ( \tilde{ \varphi }_{t_n} ( \pi ( \xi ) ) \big )  
 - \tilde{\varphi}_{t_n} ( \pi (\xi ) )  
\;=\; \varphi_{\tau } ( \tilde{x} ) - \tilde{x} \,,
\ee
that is, $\tilde{x}+\delta =\varphi_{\tau } ( \tilde{x} )$.
Therefore:
\be{bb}
 \varphi_t ( \tilde{x} + \delta ) - \varphi_t ( \tilde{x} )
\;=\; \varphi_t ( \varphi_{\tau } ( \tilde{x} )) - \varphi_t ( \tilde{x} )
\;=\; \lim_{n \rightarrow + \infty} \varphi_{t+ \tau } ( \pi \circ \varphi_{t_n} ( \xi )  )
- \varphi_t ( \pi \circ \varphi_{t_n} (\xi ) ) \,.
\ee
Now, by translation invariance, we have that:
\[
\varphi_{t+\tau }(\pi (\varphi_{t_n}(\xi ))) 
\;=\;
\varphi_{t+\tau }(\varphi_{t_n}(\xi ) - [v'\varphi_{t_n}(\xi )] v)
\;=\;
\varphi_{t+\tau }(\varphi_{t_n}(\xi )) - [v'\varphi_{t_n}(\xi )] v
\]
and similarly:
\[
\varphi_{t}(\pi (\varphi_{t_n}(\xi ))) 
\;=\;
\varphi_{t}(\varphi_{t_n}(\xi ) - [v'\varphi_{t_n}(\xi )] v)
\;=\;
\varphi_{t}(\varphi_{t_n}(\xi )) - [v'\varphi_{t_n}(\xi )] v \,
\]
so that:
\[
\varphi_{t+\tau }(\pi (\varphi_{t_n}(\xi )))
-
\varphi_{t}(\pi (\varphi_{t_n}(\xi )))
\;=\;
\varphi_{t+\tau }(\varphi_{t_n}(\xi )) - \varphi_{t}(\varphi_{t_n}(\xi )) 
\]
so, substituting into~(\ref{bb}), we have:
\be{secondremark}
\varphi_t ( \tilde{x} + \delta ) - \varphi_t ( \tilde{x} )
\;=\; \lim_{n \rightarrow + \infty}
\varphi_{t+\tau } ( \varphi_{t_n} ( \xi ) ) - \varphi_t ( \varphi_{t_n} ( \xi ) ) \;=\;
\lim_{n \rightarrow + \infty}
\varphi_{t + t_n} ( \varphi_{\tau } ( \xi ) ) - \varphi_{t + t_n } ( \xi ) \,.
\ee
Hence, (\ref{biginvariance}) follows combining (\ref{firstremark}) and
(\ref{secondremark}) (using the new sequence $\{t+t_n\}$).

Recall that $V ( \varphi_t ( \varphi_{\tau } (\xi ) ) - \varphi_t ( \xi ) )$ decreases to
its limit $\bar{V}$ as $t\rightarrow \infty $.
On the other hand, 
for any $[ \tilde{x}, \delta ] \in \Omega $, by definition of $\Omega $ we have that
$\varphi_{t_n} ( \varphi_\tau ( \xi ) ) -\varphi_{t_n} ( \xi ) \rightarrow \delta $ as $n\rightarrow \infty $.
Because of continuity of $V$, this implies that $V ( \delta ) = \bar{V}$.
Moreover, by invariance of $\Omega $, 
$V ( \varphi_t ( \tilde{x}+\delta ) - \varphi_t ( \tilde{x} ) ) = \bar{V}$, 
independently of $t$. 
Hence, application of Lemma \ref{isLyapunov} gives
$\delta \in \textrm{span} \{ v \}$ for any $[\tilde x,\delta ]\in \Omega $.
This contradicts the assumption that $\delta _0\not\in \textrm{span} \{ v \}$.
Therefore, (\ref{aa}) is true.

Projecting~(\ref{aa}) onto the $\tilde{X}$ space shows:
\[ 
\lim_{t \rightarrow + \infty}  
\tilde{\varphi}_t ( \tilde{\varphi}_\tau ( \pi (\xi ))) - \tilde{ \varphi}_t (\pi ( \xi )) \;=\; 0 \,.
\]

We next claim that every element of $\omega ( \tilde{x} )$ is an equilibrium.
Indeed, suppose that $\tilde{\varphi}_{t_n}( \pi (\xi ))\rightarrow p$;
then, for any $\tau $:
\[
\tilde \varphi_\tau (p) \;=\;
\tilde \varphi_\tau \left(\lim_{t_n \rightarrow +\infty}\tilde{\varphi}_{t_n}( \pi (\xi )) \right) \;=\;
\lim_{t_n\rightarrow + \infty} \tilde \varphi_\tau \left(\tilde{\varphi}_{t_n}( \pi (\xi )) \right) \;=\;
\lim_{t_n\rightarrow + \infty} \tilde{ \varphi}_{t_n} (\pi ( \xi )) \;=\;
p\,.
\]

Hence, the result follows by uniqueness of the equilibrium for the projected
system $\dot {\tilde x} = (I-vv')f(\tilde x)$.

\bc{competitive}
Let a system as in (\ref{thesystem}) be strongly monotone in reverse time and
enjoy the translation invariance property with respect to some vector
$v\in \interK$. Then, every solution which is bounded modulo $v$ has a
projection which converges to an equilibrium.  Moreover, there is a unique such
equilibrium. 
\ec

\bpr
The proof is entirely analogous, once Corollary \ref{isLyapunovcomp} is used
in place of Lemma \ref{isLyapunov}.
\epr

\section{An Application to Chemical Reactions}
\label{secchem}

In this section, we show how our result may be applied to conclude global
convergence to steady states, for certain chemical reactions.
A standard form for representing (well-mixed and isothermal) chemical reactions
by ordinary differential equations is:
\be{crt}
\dot {S} \;=\; \Gamma R(S),
\ee
evolving on the nonnegative orthant $\R^n_{\geq 0}$,
where $S$ is an $n$-vector specifying the concentrations of $n$ chemical
species, $\Gamma \in \R^{n\times m}$ is the \emph{stoichiometry matrix}, and
$R:\R^n_{\geq 0}\rightarrow \R^m$ is a function which provides the vector of reaction rates
for any given vector of concentrations.
We assume that $R$ is locally Lipschitz, so uniqueness of solutions holds,
and that the positive orthant $\R^n_{\geq 0}$ is invariant, and that it is
forward complete: every solution is defined for all $t\geq 0$.

To each system of the form~(\ref{crt}) and
each fixed vector $\sigma \in \R^n_{\geq 0}$, we associate the following system:
\be{f-sys}
\dot x \;=\; f_\sigma (x) \;=\; R(\sigma +\Gamma x)
\ee
evolving on the state-space
\[
X_\sigma  \;=\; \left\{x \in \R^m \,|\, \sigma +\Gamma x\geq 0 \right\} \,.
\]
The $i$th component $x_i$ of the vector $x$ is sometimes called the ``extent''
of the $i$th reaction.
We will derive conclusions about~(\ref{crt}) from the study of~(\ref{f-sys}).

Note that $X_\sigma $ is a closed set which is the closure of its interior
(it is, in fact, a polytope), and also that
$X_\sigma $ \emph{is invariant with respect to translation} by any $v\in \Ker \Gamma $,
because
$x \in X_\sigma $ 
means that $\sigma +\Gamma x\geq 0$, and therefore also $x + \lambda v\in X_\sigma $ for all $\lambda \in \R$,
because $\sigma +\Gamma (x + \lambda v) = \sigma +\Gamma x\geq 0$.

As an illustrative example, consider the following set of chemical reactions:
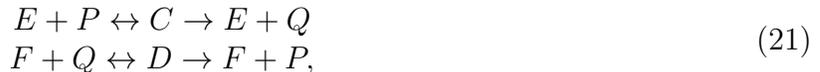
\be{singlestage}
\begin{array}{c}
E + \Su \leftrightarrow \ESu \rightarrow E + \Pd \\
F + \Pd \leftrightarrow \FPr \rightarrow F + \Su ,
\end{array}
\ee
which may be thought of as a model of the activation of a protein substrate
$\Su$ by an enzyme $E$; $\ESu$ is an intermediate complex, which dissociates
either back into the original components or into a product (activated protein)
$\Pd$ and the enzyme.
The second reaction transforms $\Pd$ back into $\Su$, and is catalyzed by
another enzyme (a phosphatase denoted by $F$).
A system of reactions of this type is sometimes called a ``futile cycle'', and
reactions of this type are ubiquitous in cell biology.
The mass-action kinetics model is obtained as follows.
Denoting concentrations with the same letters ($\Su$, etc) as the species
themselves, we introduce the species vector:
\[
S \;=\; ( \Su , \Pd , E ,  F , \ESu ,  \FPr )'
\]
and these stoichiometry matrix $\Gamma $ and vector of reaction rates $R(S)$:
\[ \Gamma = \left [ \begin{array}{cccc} -1 & 0 & 0 & 1 \\
                                         0 & 1 & -1 & 0 \\
					-1 & 1 & 0  & 0 \\
					 0 & 0 & -1 & 1 \\
					 1 & -1 & 0 & 0 \\
					 0 & 0 & 1 & -1 
\end{array}
\right ] \qquad R(S) = \left [ 
\begin{array}{c}
 k_1 E \Su - k_{-1} \ESu \\
 k_2 \ESu \\
 k_3 F \Pd - k_{-3} \FPr \\
 k_4 \FPr
\end{array} \right ] \,.
\]  
The reaction constants $k_i$, with $i=-1,1,2,3,-3,4$, are arbitrary positive
real numbers, and they quantify the speed of the different reactions.
This gives a system~(\ref{crt}).
Note that, along all solutions, one has that
\[
\Su (t)+ \Pd (t)+  \ESu (t)+   \FPr(t) \equiv  \mbox{constant}
\]
because $(1,1,0,0,1,1)\Gamma =0$.  Since the components are
nonnegative, this means that, for any solution, each of $\Su(t)$, $\Pd(t)$,
$\ESu(t)$, and $\FPr(t)$ are upper bounded by the constant
$\Su(0)+\Pd(0)+\ESu(0)+\FPr(0)$.
Similarly, we have two more independent conservation laws:
\[
E(t)+\ESu(t) \quad\mbox{and}\quad F(t)+\FPr(t)
\]
are also constant along trajectories, so also $E$ and $F$ remain bounded.
Therefore, all solutions are bounded, and hence, in particular, are defined
for all $t\geq 0$.
The system of equations~(\ref{crt}) in this example is not monotone, at
least with respect to any orthant order.
(See~\cite{monotonereactionnetworks} for more on this example, as well as
an alternative way to study it.)
We will prove, as a corollary of our main theorem, that every solution
that starts with $E(0)+\ESu(0)\not= 0$ and $F(0)+\FPr(0)\not= 0$ converges to
a steady state, which is unique with respect to the conservation relations.

\bl{lemma-R}
The system~(\ref{f-sys}) is forward complete: every solution is defined for
all $t\geq 0$ and remains in $X_\sigma $.
Furthermore, if it holds that every solution of~(\ref{crt}) is bounded, then,
for every solution $x(t)$ of~(\ref{f-sys}), $\Gamma x(t)$ is bounded. 
\els

\bpr
Pick any $x_0\in X_\sigma $, and let $S_0:=\sigma +\Gamma x_0\in \R^n_{\geq 0}$.
Consider the solution of $S(t)$ of the initial value problem
$\dot {S}= \Gamma R(S)$, $S(0)=S_0$, which is well-defined and satisfies $S(t)\geq 0$
for all $t\geq 0$.
Let, for $t\geq 0$:
\be{extent-from-reactions}
x(t)\,:=\; x_0 + \int_0^t R(S(\tau ))\,d\tau  \,.
\ee
Note that $\dot x(t)=R(S(t))$ for all $t$.
We claim that $x$ is a solution of $\dot x=f_\sigma (x)$.  Since $x(0)=x_0$ and
$x$ is defined for all $t$, uniqueness of solutions ($f_\sigma $ is locally
Lipschitz) will prove the first statement of the lemma.
To prove the claim, we first introduce the new vector function
\[
P(t)\,:=\; \sigma +\Gamma x(t)\,.
\]
Differentiating with respect to time we obtain that
$\dot P(t)=\Gamma \dot x(t)=\Gamma (R(S(t)))=\dot S(t)$ for all $t\geq 0$.  Therefore, $P-S$ is
constant.  Since $P(0)=\sigma +\Gamma x_0=S(0)$, it follows 
that $P\equiv S$.  In other words, $S$ satisfies $S(t)=\sigma +\Gamma x(t)$.
Thus,
$\dot x(t)=R(S(t))=R(\sigma +\Gamma x(t)) = f_\sigma (x(t))$, as claimed.

To prove the second statement, we simply remark that, as already proved,
for every solution $x$ of~(\ref{f-sys}), there is a solution $S$ 
of~(\ref{crt}) such that $S(t)=\sigma +\Gamma x(t)$.  Therefore,
$\Gamma x(t)=S(t)-\sigma $ is bounded if $S(t)$ is.
\epr

Note that the futile cycle example discussed earlier satisfies the assumptions
of this Lemma.  We now specialize further, imposing additional conditions also
satisfied by the example.

\bl{relate-bounded}
Suppose that the matrix $\Gamma $ has rank exactly $n-1$, its kernel spanned by
some positive unit vector $v$.
Let $x(t)$ be a solution of~(\ref{f-sys}). Then, $\Gamma x(t)$ is bounded if and
only if $\pi _vx(t)$ is bounded.
\els

\bpr
Since $\Gamma \pi _vx = \Gamma (x-(v'x)v)=\Gamma x$, one implication is clear.
Let $M$ be the restriction of $\Gamma $ to the space $v^{ \perp}$ orthogonal to
the vector $v$, i.e.\ the image of $\pi _v$.
As $\Gamma \pi _vx =\Gamma x$, the images of $\Gamma $ and $M$ are the same.
The map $M$ is one-to-one: suppose that $x\in v^{ \perp}$
is so that if $Mx=0$.   Then, $\Gamma x=0$, so $x$ is in the kernel of $\Gamma $,
i.e., it is also in the span of $v$.  Thus, $x=0$.
Let $M^{-1}$ be the inverse of $M$, mapping the image of $\Gamma $ into $v^{\perp}$.
Thus, if a trajectory is such that $\Gamma x(t)$ is bounded, then also
\[
M^{-1}\Gamma x(t) = M^{-1}\Gamma \pi _vx(t)=  M^{-1}M\pi _vx(t) = \pi _vx(t)
\]
is bounded. 
\epr

Observe that the spaces $X_\sigma $ are translation invariant with respect
to any $v$ as in the statement of this Lemma.

\bc{chem-cor}
Suppose that:
\ben
\item
the matrix $\Gamma $ has rank $n-1$, with kernel spanned by
some positive unit vector;
\item
every solution of~(\ref{crt}) is bounded;
\item
$\sigma \in \R^n_{\geq 0}$ is so that
the system $\dot x=f_\sigma (x)$ is strongly monotone.
\een
Then, there is a $\zeta =\zeta _\sigma \in \R^n_{\geq 0}$ with the following property:
for each $\rho \in \R^n_{\geq 0}$ such that $\rho -\sigma \in \mbox{Image}(G)$,
the solution $S$ of~(\ref{crt}) with $S(0)=\rho $ satisfies
$S(t)\rightarrow \zeta $ as $t\rightarrow \infty $.
\ecs

\bpr
We let the kernel of $\Gamma $ be spanned by the positive unit vector $v$.
By Lemmas~\ref{lemma-R} and~\ref{relate-bounded}, $\pi _vx(t)$ is bounded, for
every solution of~(\ref{f-sys}).
By Theorem~\ref{main}, there is a unique equilibrium 
$\xi $ of the projected system $\dot {\widetilde x} = (I-vv')f(\widetilde x)$
so that every solution $x$ of $\dot x=R(\sigma +\Gamma x)$ is such that
$\pi _v(x(t))\rightarrow \xi $ as $t\rightarrow \infty $.
We next show that $\zeta =\sigma +\Gamma \xi $ satisfies the requirements.

Pick $\rho \in \R^n_{\geq 0}$ so that $\rho -\sigma =\Gamma a$, $a\in \R^m$, and let
$S$ be the solution of $\dot S=\Gamma R(S)$ with initial condition $S(0)=\rho $.
Arguing as in the proof of Lemma~\ref{lemma-R}, we have that
$S(t)=\rho +\Gamma x(t)$, where
$\dot x=R(\rho +\Gamma x)$, $x(0)=0$.

Introduce the function $z(t)=x(t)+a$.
Then, $\dot z=\dot x+0=R(\rho +\Gamma x)=R(\sigma +\Gamma z)$, with $z(0)=a$.
Since $\sigma +\Gamma z(0)=\sigma +\Gamma a=\rho \geq 0$, it follows that $z(0)\in X_\sigma $, and therefore
$z(t)$ is a solution of $\dot x=R(\sigma +\Gamma x)$ on $X_\sigma $.
Therefore, $\pi _vz(t)\rightarrow \xi $.
As $x(t)=z(t)-a$, this means that $\pi _vx(t)\rightarrow \xi -\pi _va$.
Since for every vector $x$ it holds that $\Gamma \pi _vx=\Gamma x$, applying $\Gamma $
to the above gives
\[
\Gamma x(t) = \Gamma \pi _vx(t) \rightarrow  \Gamma \xi  -\Gamma a\,.
\]
Therefore, 
$S(t)=\rho +\Gamma x(t)\rightarrow \rho +\Gamma \xi -\Gamma a=\sigma +\Gamma \xi =\zeta $ as $t\rightarrow \infty $.
\epr

In the futile cycle example, we may take $v=(1/4,1/4,1/4,1/4)'$, and consider
the following set:
\[
\Sigma  = \{\sigma = (\Su,\Pd,E,F,\ESu,\FPr)\in \R^n_{\geq 0} \st
E+\ESu>0, F +\FPr>0 \}\,.
\]
The system $\dot x=f_\sigma (x)$ is strongly monotone for $\sigma \in \Sigma $.  To see this,
we compute the Jacobian of $R(\sigma  + \Gamma x(t))$ with respect to $x$:
\[
\pmatrix{* &     * & 0 & k_1E\cr
         * &     * & 0 & 0   \cr
         0 &  k_3F & * & *   \cr
         0 &     0 & * & *   }
\]
where the stars represent strictly positive elements when in the off-diagonals
(and strictly negative when on the diagonals), and where $E,F$ are the
$E$ and $F$ coordinates of $\sigma  + \Gamma x$, or, more explicitly:
\[
\pmatrix{* &     * & 0 & k_1(\sigma _3+(x_2-x_1))\cr
         * &     * & 0 & 0   \cr
         0 &  k_3(\sigma _4+(x_4-x_3)) & * & *   \cr
         0 &     0 & * & *   } \,.
\]
Thus, the system is cooperative (i.e., monotone with respect to the main
orthant).  
It is strongly monotone if this matrix is irreducible almost everywhere along
trajectories (see e.g.~\cite{Hirsch-Smith}, Section~3.2), which
amounts, because $f_\sigma $ is a real-analytic function, to asking that 
$\sigma _3+x_2-x_1\not\equiv 0$ and $\sigma _4+x_4-x_3\not\equiv 0$ along any solution.
Let us prove now that this is the case, assuming that $\sigma \in \Sigma $,
that is, that $\sigma _3+\sigma _5\not= 0$ and $\sigma _4+\sigma _6\not= 0$.
Suppose that $\sigma _3+x_2-x_1 \equiv  0$, so that $\dot x_1-\dot x_2\equiv 0$ and
$x_1-x_2\equiv \sigma _3$.
The equations for~(\ref{f-sys}) give:
\[
\dot x_1 - \dot x_2
\;=\;
k_1(\sigma _3+x_2-x_1)(\sigma _1+x_4-x_1) - (k_{-1}+k_2)(\sigma _5x_1-x_2) \,,
\]
so:
\[
0 \equiv  - (k_{-1}+k_2)(\sigma _3+\sigma _5)
\]
which contradicts $\sigma _3+\sigma _5\not= 0$.
Similarly for $\sigma _4+x_4-x_3 \equiv  0$. 
So the system is indeed strongly monotone.

We conclude that every solution of our example with an initial condition
in the set $\Sigma $ converges to an equilibrium.  Moreover, there is a unique such
equilibrium in each stoichiometry class $\sigma +\mbox{Image}(\Gamma )$.

When initial conditions do not belong to $\Sigma $, one has a standard enzymatic
Michaelis-Menten type of reaction, and the same conclusion holds.
This is very easy to show.
(Indeed, take for instance the case when $E(0)=C(0)=0$.
As $\dot P=k_4D$, $P(t)$ is nondecreasing, so (since it is upper bounded)
we know that $P$ converges.
Consider the function $y=Q+D$.  Since $P+y$ is constant, $y$ converges, too.
Since $\dot y$ has a bounded derivative (it can be expressed in terms
of bounded variables), and its integral is convergent, it follows
(``Barbalat's lemma'') that $\dot y=-k_4D$ converges to zero, so $D$ must converge
and therefore, again using that $P+Q+D$ is constant, $Q$ converges as well.
Finally, since $D+F$ is constant, $F$ converges, too.)

\section{Remarks on Duality and Possible Extensions}

As pointed out in the introduction, our main result stated in Theorem
\ref{main} can be seen as a dual to Mierczynski's global convergence theorem
for strongly cooperative systems with a positive first integral, published in
\cite{mier}. 
We discuss this informally in this section.
Strictly speaking, duality of \ref{main} only holds provided that we
consider the following special case of Mierczinski's Theorem:
{\em
Consider a system of ordinary differential equations in $\R^n_{+}$, defined by
a $\mathcal{C}^1$ vector field $f: \R^n_+ \rightarrow \R^n$, such that:
$f(0)=0$,
$\frac{\partial f_i}{\partial x_j} > 0$ for all $i \neq j$,
and:
\be{dual-eqn}
\mbox{there exists a vector $c \in (\R_+)^n$ such that $c' f(x) = 0$ 
for all $x\in (\R_+)^n$.}\ee
Then, every solution is bounded and converges to an equilibrium.}
This is a special case of Mierczinski result, which had already appeared in
several previous publications, in that \emph{linear} positive first integrals
are considered; namely the quantity $c'x$ is preserved along solutions of the
system.  For simplicity, we actually strengthened 
one of the original assumptions by asking that $\frac{\partial f_i}{\partial x_j} > 0$ for
all $i \neq j$, rather than a strict monotonicity condition with respect to
all off-diagonal entries as the Theorem is stated in \cite{mier}.

The duality with Theorem \ref{main} is evident if we express the conditions in
terms of the Jacobian of the vector field.
A linear positive first integral amounts to having a constant
left-eigenvector relative to the dominant zero eigenvalue for the Jacobian
matrix $D f(x)$; in particular, $c^{'} Df(x) = 0$ for all $x$ in the
state-space.
On the other hand, translation invariance by a positive vector $v$ (over a
given state-space $X$) can be stated in terms of the Jacobian matrix by asking
that $D f(x) v =0$ for all $x \in X$; i.e., the existence of a constant
right-eigenvector relative to the dominant zero eigenvalue of the Jacobian
matrix $Df(x)$. 
As a further remark, we note that our main result does not need the strict
monotonicity condition as stated above, but only asks for strong-monotonicity
of the resulting flow (this is in fact weaker than assuming strictly positive
off-diagonal entries of the Jacobian; for instance a much tighter sufficient
condition for strong monotonicity of the flow, can be formulated by asking
that the Jacobian matrix have \emph{non-negative} off-diagonal entries and be
\emph{irreducible}).

A more general version of Mierczinski Theorem than stated above does not
assume linearity of the first integral.
In particular, assumption~(\ref{dual-eqn} is replaced by the existence of a
$\mathcal{C}^1$ function  $H(x):\R^n_+ \rightarrow \R $, such that $D H(x) \cdot f(x) =0$
and $DH(x) \in \R^n_+$ for all $x \in \R^n_+$. 
This condition does not allow an elegant interpretation in terms of Jacobians
of $Df(x)$, but nevertheless, one may state a nonlinear dual of
the Theorem would provided that we understand translation invariance in the
following more general sense.
Let us say
that a flow is \emph{invariant with respect to translation by a strictly
  increasing flow $\tilde{\varphi}$}
if for all $t_1$,$t_2$ in $\R$ the following holds:
\[ 
\varphi_{t_1} ( \tilde{ \varphi }_{t_2} (x_0) ) = \tilde{ \varphi}_{t_2} ( \varphi_{t_1} (x_0) ),
\]
and moreover, for each $x_1, x_2 \in X$ there exists $t \in \R$ so that $x_2
\succeq \tilde{\varphi}_t (x_1) $.
This property generalizes our previous concept:
translation invariance with respect to a constant vector $v$ is exactly the
property of invariance with respect to translation by the increasing flow
$\tilde{\varphi}$ induced by the system of differential equations $\dot {x} = v$.
Invariance with respect to non-trivial general flows as in this definition is
not easy to check in concrete examples, however, at least in principle, an
infinitesimal characterization of the property is as follows.
Let $f(x): X \rightarrow \R^n$ and $v(x): X \rightarrow \R^n$ be $\mathcal{C}^1$
vector-fields. The flow induced by the system $\dot {x} = f(x)$ commutes with
respect to the strictly increasing 
flow induced by $\dot {x} = v(x)$ if and only if:
\[ 
D f(x) v(x) = Dv(x) f(x) \,.
\]
Moreover, if there exists a compact set $P \subset \interK$ so that $v(x) \in P$
for all $x \in X$, then the flow induced by $v(x)$ is strictly increasing
(meaning that its solutions are such with respect to $t$) and for any $x_1$
and $x_2$ in $X$  there exists $t \in \R$ so that $x_2\succeq \tilde{\varphi}_t(x_1)$.

Accordingly we have to redefine the notion of boundedness modulo translation by $\tilde{\varphi}$ by 
asking that solutions are bounded if there exists $M>0$ such that for all $x_0 \in X$ and all $t \in \R$,
 there exists $\tau $ with the property that $| \tilde{\varphi}_{\tau } ( \varphi_t (x_0)) | \leq M$.
While this definition is rather natural, there is not, however, a natural
counter-part to the space  $\tilde{X} = X \cap v^{\perp}$.
Hence, we may as well let $\tilde{X}$ be defined
 as a quotient space of $X / \sim$ under the equivalence relation
$x_1 \sim x_2$ if and only if
 $\tilde{ \varphi}_t (x_1) = x_2$ for some $t \in \R$. 
This definition of $\tilde{X}$ and the commutativity
 of $\tilde{\varphi}$ and $\varphi$ allow us to define 
 a flow on equivalence classes of $[x]$ of $\tilde{X}$ in the natural way:
 $\varphi_t ( [x] ) := [ \varphi_t (x) ]$. 
Boundedness of a solution in the space $\tilde{X}$ is equivalent
 to boundedness modulo translation given above. 
Our main result would then be translated into the following statement
in the  current set-up:
\emph{
Consider a forward complete, strongly monotone nonlinear system
(\ref{thesystem}) with translation invariance 
with respect to a strictly increasing flow. 
Then, every solution which is bounded is such that $\varphi_t ( [x] )$ admits a
limit as $t \rightarrow + \infty $.}


\newpage
\begin{thebibliography}{99}


\itemsep0pt

\newcommand{\book}[1]{#1,}
\newcommand{\inbook}[1]{in: #1,}
\newcommand{\journal}[1]{#1}
\newcommand{\arttitle}[1]{#1,}
\newcommand{\jvol}[1]{#1}
\newcommand{\jyear}[1]{(#1)}
\newcommand{\pp}[1]{#1}
\newcommand{\AND}{}
\newcommand{\auth}[2]{#1, #2,}

\bibitem{LSUmonotone}
Angeli, D.,
Sontag, E.D.,
\arttitle{Interconnections of monotone systems with steady-state characteristics}
\inbook{Optimal Control, Stabilization, and Nonsmooth Analysis}
de Queiroz, M., M. Malisoff, and P. Wolenski (Eds.),
(Springer-Verlag, Heidelberg, 2004),
\pp{135--154}.


\bibitem{monotonereactionnetworks} 
Angeli, D., P. De Leenheer, E.D. Sontag, 
\arttitle{Monotone chemical reaction networks}
in preparation

\bibitem{craciun}
Craiciun, G.,
Feinberg, M.,
\arttitle{Multiple equilibria in complex chemical reaction networks: I. The injectivity property}
\journal{SIAM Journal on  Applied Mathematics}
\jvol{65}
\jyear{2005}
\pp{1526--1546}.

\bibitem{Dancer}
Dancer, E.N.,
\arttitle{Some remarks on a boundedness assumption for monotone dynamical systems}
\journal{Proc. of the AMS}
\jvol{126}\jyear{1998}
\pp{801--807}.

\bibitem{feinberg}
Feinberg, M.,
\arttitle{Chemical reaction network structure and the stabiliy of complex isothermal reactors - I. The deficiency zero and deficiency one theorems}
\journal{Chemical Engr. Sci.}
\jvol{42}\jyear{1987}
\pp{2229--2268}.

\bibitem{gouze}
Gouz\'e, J.-L.,
\arttitle{Global behaviour of Lotka-Volterra systems}
\journal{Mathematical Biosciences}
\jvol{113}\jyear{1993}
\pp{231--243}.

\bibitem{Hirsch}
Hirsch, M.,
\arttitle{Differential equations and convergence almost everywhere in strongly monotone flows}
\journal{Contemporary Mathematics}
\jvol{17}\jyear{1983}
\pp{267--285}.

\bibitem{Hirsch2}
Hirsch, M.,
\arttitle{Systems of differential equations that are competitive or cooperative II: Convergence almost everywhere}
\journal{SIAM J. Mathematical Analysis}
\jvol{16}\jyear{1985}
\pp{423--439}.

\bibitem{hirsch} 
Hirsch, M.,
\arttitle{Stability and convergence in strongly monotone dynamical systems}
\journal{J. Reine Angew. Math.}
\jvol{383}\jyear{1988}
\pp{1--53}.

\bibitem{Hirsch-Smith}
Hirsch, M.,
H.L. Smith,
\arttitle{Monotone dynamical systems}
\inbook{Handbook of Differential Equations, Ordinary Differential
Equations (second volume)}
Elsevier, 2005, to appear.

\bibitem{JiFa}
Jiang, J.F.,
\arttitle{On the global stability of cooperative systems}
\journal{Bulletin of the London Math Soc}
\jvol{6}\jyear{1994}
\pp{455--458}.

\bibitem{kunze} 
\auth{Kunze}{H.}
\AND
\auth{Siegel}{D.}
\arttitle{Monotonicity properties of chemical reactions with a single initial bimolecular step}
\journal{J.\ Math.\ Chem.}
\jvol{31}\jyear{2002}
\pp{339--344}.

\bibitem{mier} 
Mierczynski, J.,
\arttitle{Strictly cooperative systems with a first integral}
\journal{SIAM J. Math. Anal.}
\jvol{18}\jyear{1987}
\pp{642--646}.

\bibitem{Smillie}
J.\ Smillie,
\arttitle{petitive and cooperative tridiagonal systems of differential equations}
\journal{SIAM J. Math. Anal.}, 
\jvol{15}\jyear{1984}
\pp{530--534}.

\bibitem{Smith}
H.L.\ Smith, 
\book{Monotone dynamical systems: An introduction to the theory of competitive and cooperative systems, Mathematical Surveys and Monographs, vol. 41}
AMS, Providence, RI, 1995.

\bibitem{sysbio}
Sontag, E.D.
\arttitle{Some new directions in control theory inspired by systems biology}
\journal{Systems Biology}
\jvol{1}\jyear{2004}
\pp{9--18}.

\bibitem{volpert}
\auth{Vol'pert}{A.}
\AND
\auth{Hudjaev}{S.}
\book{Analysis in Classes of Discontinuous Functions and Equations of Mathematical Physics}
Marinus Nijhoff, Dordrecht, 1985.


\end{thebibliography}
\end{document}